\documentclass{amsart}
\usepackage[headings]{fullpage}
\usepackage{ amssymb, amsfonts, amsmath, color}

\newtheorem{thm}{Theorem}

\newtheorem{prop}[thm]{Proposition}
\theoremstyle{definition}

\theoremstyle{remark}
\newtheorem*{rem}{Remark}

\newcommand{\LM}{{\mathcal{LM}}}
\newcommand{\g}{{\mathfrak{g}}}

\title{{A comment on the integration of Leibniz algebras}}
\author{J.Mostovoy} 

\begin{document}

\begin{abstract}
In this note we point out that the definition of the universal enveloping dialgebra for a Leibniz algebra is consistent with the interpretation of a Leibniz algebra as a generalization not of a Lie algebra, but of the adjoint representation of a Lie algebra. From this point of view, the formal integration problem of Leibniz algebras is, essentially, trivial.
\end{abstract}

\maketitle

A Leibniz algebra is a vector space with a bilinear bracket satisfying the Leibniz identity
$$[[x,y],z]=[[x,z],y]+[x,[y,z]].$$ 
A Leibniz algebra whose bracket is antisymmetric is a Lie algebra. 
J.-L.\ Loday who defined Leibniz algebras suggested that they arise as tangent structures to some hypothetic objects he called ``coquecigrues'' in the same way as the Lie algebra structure arises
on the tangent space to a Lie group at the unit element. 
The hunt for coquecigrues has resulted in various trophies such as \cite{KW}. However, the unsatisfying feature of the integration method in \cite{KW} is that in the case when the Leibniz algebra is a Lie algebra it does not produce a Lie group. 

The purpose of this note is point out that this feature is consistent with the definition of the universal enveloping dialgebra of a Leibniz algebra. Namely, we show that the usual procedure of the formal Lie integration via distribution algebras, as described in \cite{Serre}, can be applied in the context of Leibniz algebras. For Leibniz algebras with an antisymmetric bracket, that is, for Lie algebras, this integration method produces not the corresponding Lie group but, rather, the group together with its adjoint representation.
It will be more convenient for us to replace Leibniz algebras by a wider class of objects, namely, Lie algebras in the Loday-Pirashvili category of linear maps (2-vector spaces). We shall see that such a Lie algebra integrates to a triple that consists of a Lie group $G$, its right representation $\rho$ and a morphism of $\rho$ to the right adjoint representation $\mathrm{Ad}^{-1}_G$.

In a way, the output of this integration procedure is trivial. A Leibniz algebra can be thought of as a right $\g$-module over the right adjoint representation of $\g$ for some Lie algebra $\g$ and this obviously integrates to a right representation of the corresponding Lie group $G$ with a morphism to the right adjoint representation of $G$. The main result of this note says that this kind of integration is the natural consequence of considering dialgebras as the category of universal enveloping objects for Leibniz algebras. 

It is worth pointing out that we do not claim here that Leibniz algebras cannot be integrated so as to produce Lie groups in the special case of Lie algebras. Indeed, a {\em local} integration procedure of this type was described by S.\ Covez in \cite{C}, who developed the idea of M.\ Kinyon \cite{K}. It is shown in \cite{C} that any Leibniz algebra can be locally integrated to a local augmented Lie rack; this integration method applied to a Lie algebra produces a local Lie group. The algebraic structure on the distributions  supported at the unit of a local augmented Lie racks should produce another version of the universal enveloping algebra for Leibniz algebras. 

All vector spaces and algebras will be assumed to be defined over a field of characteristic zero.

\section{Lie algebras in the Loday-Pirashvili category}
\subsection{The Loday-Pirashvili category}
The Loday-Pirashvili category $\LM$  (or the category of 2-vector spaces), defined in \cite{LP}, has as objects linear maps $(V\to W)$ between vector spaces. A morphism from $(V\to W)$ to $(V'\to W')$ 
is a commutative diagram
$$
\begin{array}{ccc}
V&\to&V'\\
\downarrow&&\downarrow\\
W&\to&W'.
\end{array}
$$
There is a tensor product defined as
 $$(V\xrightarrow{\delta} W)\otimes(V'\xrightarrow{\delta'} W')=((V\otimes W' + W\otimes V')\xrightarrow{\delta\otimes \text{id}+\text{id}\otimes \delta'} W\otimes W'),$$ 
 where, for the sake of readability, we write $+$ for the direct sum. If objects of $\LM$ are thought of as chain complexes concentrated in degrees 1 and 0, this tensor product is just the product of chain complexes with the degree 2 part thrown away. We shall refer to the vector spaces $V$ and $W$ as the degree 1 and degree 0 parts of $(V\to W)$ respectively.

The interchange automorphism $\tau$ defined on $(V\to W)^{\otimes 2}$ acts by
$$\tau( x\otimes y' + y\otimes x')= x'\otimes y + y'\otimes x$$
in degree 1 and  by $\tau(x\otimes x')=x'\otimes x$ in degree 0. The interchange automorphism gives an action of the symmetric group $\Sigma_n$ on the $n$th tensor power  of an object in $\LM$. This allows to define symmetric and exterior powers:  the symmetric power is universal for the morphisms from the tensor power of that are invariant under this action, and the exterior power  for the morphisms that change by the sign representation of $\Sigma_n$.

In particular, $$S^k(V\to W)=(S^{k-1}W\otimes V \to S^k W),$$ and 
$$\Lambda^k(V\to W)=(\Lambda^{k-1}W\otimes V \to \Lambda^k W),$$ 
see \cite{LP} for details. It is easy to see that, similarly to the isomorphism $S(W+W')=S(W)\otimes S(W')$ there is an isomorphism
$$S\left((V+V')\xrightarrow{\delta+\delta'}(W+W')\right)=S(V\xrightarrow{\delta} W)\otimes S(V'\xrightarrow{\delta'} W')$$
for any pair of objects $(V\xrightarrow{\delta} W)$ and $(V'\xrightarrow{\delta'} W')$ in $\LM$.

\subsection{Lie algebras in $\LM$ and their universal enveloping algebras}
A Lie algebra in $\LM$ is a linear map $\mu:(V\to W)^{\otimes 2}\to (V\to W)$ which is antisymmetric (that is, which factors through $\Lambda^2(V\to W))$ and satisfies the Jacobi identity
$$\mu(1\otimes\mu)-\mu(\mu\otimes 1)+\mu(\mu\otimes 1)(1\otimes\tau)=0.$$
In general, a Lie algebra is $\LM$ is an object $(M\xrightarrow{\delta}\g)$ with $\g$ a Lie algebra, $M$ a right $\g$-module and $\delta$ a $\g$-equivariant map. A Leibniz algebra $\g$ gives rise to a Lie algebra $(\g\to\g_{Lie})$ in $\LM$. Conversely, given a Lie algebra $(M\xrightarrow{\delta}\g)$ in $\LM$ one can define the Leibniz algebra bracket on $M$ by 
$$[x,y]=[x,\delta y],$$
where the bracket on the right-hand side denotes the right action action of $\g$ on $M$, see \cite{LP}.

Similarly, one can speak of associative algebras, coalgebras in $\LM$ and so on. The universal enveloping algebra of a Lie algebra $(M\to\g)$ in $\LM$  is the $U(\g)$-bimodule $U(\g)\otimes M$, such that for all $g\in \g$, $h\in U(\g)$ and $m\in M$
$$g\cdot h\otimes m= gh\otimes m$$
and
$$h\otimes m \cdot g= hg\otimes m+ h\otimes [m,g];$$
here $[m,g]$ denotes the right action of $g$ on $M$.

The degree 1 part of an algebra $(B\xrightarrow{\delta} A)$ in $\LM$ carries the structure of a dialgebra (see \cite{L2})  given by 
$$ x \vdash y = \delta x\cdot y \quad \text{and} \quad x \dashv y = x \cdot \delta y $$
for $x,y\in B$. The universal enveloping dialgebra of a Leibniz algebra $\g$ is the dialgebra structure on  $U(\g_{ie})\otimes \g$ coming from the universal enveloping algebra of the Lie algebra $(\g\to\g_{Lie})$ in $\LM$, see \cite{L2, G}.

There are two kinds of trivial examples of Lie algebras in $\LM$. Any Lie algebra $\g$ gives rise to the Lie algebra $(0\to \g)$ in $\LM$, and any vector space $V$ gives rise to the Lie algebra $(V\to 0)$. The corresponding universal enveloping algebras are $(0\to U(\g))$ and  $(V\to0)$, respectively.

\subsection{The identity map of a Lie algebra and the bimodule of 1-currents}
An important and interesting class of examples of Lie algebras in $\LM$ consists of the identity maps $(\g\to\g)$, where $\g$  is a Lie algebra. These are precisely the Lie algebras considered as Leibniz algebras. In this case, the universal enveloping algebra can be interpreted in terms of the bimodule of 1-currents on a local analytic Lie group $G$ whose algebra is $\g$. 

Recall (\cite{Serre}) that the algebra $D_0(G)$ of distributions on $G$ supported at the unit element is naturally isomorphic to $U(\g)$. As a vector space it is isomorphic to the symmetric algebra $S(\g)$ by the Poincar\'e-Birkhoff-Witt Theorem. Indeed, each distribution in $D_0(G)$ can be though of a differential operator applied to the Dirac's delta function. 

Denote by $D_1(G)$ the space of 1-currents (that is, linear functionals on 1-forms) on $G$ 
supported at the unit. By definition, each 1-current in $D_1(G)$ can be written as a sum\footnote{We shall use the Einstein summation convention.} $\alpha_i\partial_i$ where $\alpha_i$ are distributions supported at the unit and $\partial_i$ is the coordinate in $\g$ dual to $dx_i$. Applied to a 1-form $f_j dx_j$ the current   $\alpha_i\partial_i$ gives the sum $\alpha_i(f_i)$.
As a vector space, $D_1(G)$ is isomorphic to $S(\g)\otimes\g$. 


\medskip

There is a linear map $$\delta: D_1(G)\to D_0(G),$$
given by $$(\delta\alpha) (f)=\alpha (df)$$
for $\mu\in D_1(G)$ and $f$ a 
function on $G$.

\medskip

The product $\mu:G\times G\to G$ and the diagonal $\Delta: G\to G\times G$ induce maps of the spaces of 
differential forms $$\mu^*: \Omega^1(G)\to\Omega^1(G\times G)=\Omega^1(G)\,\widehat{\otimes}\, \Omega^0(G) + \Omega^0(G)\ \widehat{\otimes}\ \Omega^1(G)$$
and 
$$\Delta^*:\Omega^1(G\times G)=\Omega^1(G)\ \widehat{\otimes}\ \Omega^0(G) + \Omega^0(G)\ \widehat{\otimes}\ \Omega^1(G)\to \Omega^1(G).$$ Dually, there are maps
$$\mu':D_1(G)\otimes D_0(G) + D_0(G)\otimes D_1(G)\to D_1(G)$$
and
$$\Delta': D_1(G)\to D_0(G)\otimes D_1(G) + D_1(G)\otimes D_0(G),$$
which give the map $\delta: D_1(G)\to D_0(G)$ the structure of a bialgebra in $\LM$.

The primitive elements in degree 1 are the currents of the form $\alpha_i\partial_i$ where each $\alpha_i$ is a constant and $\partial_i$ is dual to $dx_i$; in degree 0 the primitives are of the form  $\alpha_i\partial_i$ where  each $\alpha_i$ is a constant and $\partial_i$ is the derivative of the Dirac's delta along the $i$th coordinate axis, evaluated at the origin. The map $\delta$ identifies both primitive subspaces. By the Milnor-Moore Theorem of \cite{LP}, we have
\begin{prop}
The map $(D_1(G)\xrightarrow{\delta} D_0(G))$ is naturally isomorphic to the universal enveloping algebra of $(\g\xrightarrow{\mathrm{id}}\g)$. 
\end{prop}
In particular, the universal enveloping dialgebra of a Lie algebra has a fundamentally different geometric meaning from its universal enveloping algebra: it consists of 1-currents and not of distributions.


\section{Integration}
\subsection{Formal integration}
Recall that a formal group on a vector space $V$ is a linear map
$$F: S(V+V)\to V,$$
which is associative and unital. The fact that $F$ is unital means
$$F|_{1\otimes S(V)}=1\otimes \pi_V\quad \text{and}\quad F|_{S(V) \otimes 1}= \pi_V \otimes 1,$$ where $\pi_V:S(V)\to V$ is the projection onto the degree 1 subspace. Associativity means that the extension of $F$ to a coalgebra morphism $$F':S(V+V)=S(V)\otimes S(V)\to S(V)$$
is an associative product. (Any linear map $\theta: S(V)\to V'$ can be extended to a unique coalgebra morphism $\theta':S(V)\to S(V')$. The extension is given explicitly by the formula
$$\theta'(\mu)=\sum_{n=0}^{\infty}\frac{1}{n!}\,\theta(\mu_{(1)})\ldots \theta(\mu_{(n)})=\epsilon(\mu)1+\theta(\mu)+\ldots,$$
where $\epsilon$ is the counit and Sweedler's notation is used.)

The map $F$ is interpreted as follows in terms of an $n$-tuple  $(f_k(x_i,y_j))$ of power series which represent the formal group ($n=\dim V$).
Choose a basis in $V$ and let  $x_i$ and $y_j$ be the coordinates in the first and the second copies of $V$ respectively, and $F_k$ the components of the map $F$. Then $F_k$ sends a monomial in $x_i$ and $y_j$ to  its coefficient in $f_k$.

Given a Lie algebra $\g$, a formal group which  integrates $\g$ can be obtained via the Campbell-Baker-Hausdorff formula. Alternatively, consider the map
$$U(\g)\otimes U(\g)\to U(\g)\xrightarrow{\text{Prim}} \g,$$
where the first arrow is the product in the universal enveloping algebra and the second arrow is the projection onto the primitive subspace. Identifying $U(\g)$ with $S(\g)$ via the Poincar\'e-Birkhoff-Witt Theorem, we obtain a formal group which integrates $\g$.

This approach to formal integration can be applied in tensor categories other than vector spaces; in particular, in the Loday-Pirashvili category.

Define a formal group in the Loday-Pirashvili category to be an object $(V\xrightarrow{\delta} W)$ of $\LM$ together with a linear map
$$S\left((V+V)\xrightarrow{\delta+\delta}(W+W)\right)\xrightarrow{G} \left(V\xrightarrow{\delta} W\right),$$
whose extension to a coalgebra  morphism in $\LM$
$$S(V\xrightarrow{\delta}W)\otimes S(V\xrightarrow{\delta}W) \to  S(V\xrightarrow{\delta}W) $$
is an algebra in $\LM$. Here we use the fact that, as in the case of usual coalgebras, any linear morphism $\theta: S(V\to W)\to (V'\to W')$ can be extended to a unique coalgebra morphism.  In degree 1 this extension is given by 
$$\theta'_1(\mu\otimes v)=
\sum_{n=0}^{\infty}\frac{1}{n!}\,\theta_0(\mu_{(1)})\ldots \theta_0(\mu_{(n)})\otimes \theta_1\left(\mu_{(n+1)}\otimes v\right),$$
 where $\theta_i$ is the morphism between the degree $i$ components.)

With this definition of a formal group the integration problem is, essentially, trivial. 
Given a Lie algebra $(M\to\g)$ in $\LM$, compose the product in $U(M\to\g)$ with the projection to the primitive subspace. Identifying $U(M\to\g)$ with $S(M\to\g)$ we get a diagram
$$
\begin{array}{ccc} 
\medskip S(\g)\otimes M\otimes S(\g) + S(\g)\otimes S(\g)\otimes M& \xrightarrow{G^1+G^2} &M\\
\medskip \downarrow&&\downarrow\\
S(\g)\otimes S(\g)& \xrightarrow{F} &\g
\end{array}
$$ 
which ``integrates'' the Lie algebra $(M,\g)$ in $\LM$. 
Conversely, given a formal group, its extension to a coalgebra morphism is a bialgebra whose primitive subspace is a Lie algebra in $\LM$. We have

\begin{prop}
The functor that assigns to a Lie algebra $(M\to\g)$ in $\LM$ the primitive part of the product in $U(M\to\g)$ is an equivalence of the categories of Lie algebras in $\LM$ and of formal groups in $\LM$. 
\end{prop}

Clearly, a Lie algebra $(0\to\g)$ integrates to a usual formal group and a vector space $(V\to 0)$ integrates to itself.  For a general Lie algebra $(M\xrightarrow{}\g)$ in $\LM$, consider first the primitive part of the left action of $U(\g)$  on $U(\g)\otimes M$:
$$\alpha\otimes (\beta\otimes v) \to \alpha\beta\otimes v \xrightarrow{\text{Prim}} \epsilon(\alpha\beta) v,$$
where $\alpha,\beta$ are in $U(\g)$, $v\in M$ and $\epsilon$ is the counit in $U(\g)$. This means that if this map is thought of as an $m$-tuple, where $m=\dim M$, of formal power series $g^1(x,y,v)$ with $x,y\in\g$ and $v\in M$, it is of a very simple form
$$g^1(x,y,v)=v,$$
which neither depends on the $\g$-module structure of $M$, nor on the Lie algebra structure of $\g$.  
As for the primitive part of the right action of $U(\g)$  on $U(\g)\otimes M$ we have
$$(\alpha\otimes v)\otimes 1 \to \alpha\otimes v \xrightarrow{\text{Prim}}  \epsilon(\alpha) v,$$
and for $b\in \g$
$$(\alpha\otimes v)\otimes b\to\alpha b\otimes v + \alpha\otimes [v,b]\xrightarrow{\text{Prim}}  \epsilon(\alpha) [v,b].$$
In particular, we see that the primitive part of the right action does not depend on the non-zero degree terms in $\alpha$. As an $m$-tuple of  formal power series, this action can be written as a linear in $v$ function $g^2(x,v,y)=g^2(v,y)$. 

There are two conditions on the function $g^2(v,y)$. Let $f(x,y)$ be the $n$-tuple ($n=\dim \g$) of formal power series representing the formal group $F$ on $\g$. Then the first condition is associativity:
$$\g^{2}(g^{2}(v,y),z)=\g^{2}(v,f(y,z)),$$
and the second is the relation to the function $f$: $\delta g^2(v,y)$ should coincide with the linear in $\delta v$ terms of $f(\delta v, y)$; here $\delta$ is the map $M\to\g$.

For example, for the Lie algebra $(\g\xrightarrow{\text{id}}\g)$ the function $g^{2}(v,y)$ is simply the linear in $v$ part of $f(v,y)$.

\subsection{A global interpretation} The formal group in $\LM$ that integrates the identity map $(\g\to\g)$ can be thought of as an infinitesimal version of the product on the tangent bundle of a Lie group $G$:
$$TG\times TG\to TG.$$
This product is, in fact, a pair of actions (right and left) of $G$ on $TG$. If the left action is used to trivialize $TG$, the right action becomes the right adjoint action $\mathrm{Ad}^{-1}$ of $G$ on $\g$. 

The global version of an arbitrary formal group in $\LM$ is a vector bundle $\xi$  over $G$ with the fibre $M$, together with an ``anchor map'' $p:\xi\to TG$ which commutes with the projections, and a pair of actions, right and left, of $G$ on $\xi$ which are carried by $p$ to the actions of $G$ on $TG$. The bundle $\xi$ can be trivialized by means of the left action and we see that a ``Lie group in $\LM$'' is simply a commuting triangle
\begin{equation}\label{LPG}
\begin{array}{ccc}
G&\xrightarrow{\mathrm{\rho}}& GL(M)\\
\Vert& & \downarrow\\
G&\xrightarrow{\mathrm{Ad^{-1}}}&GL(\g)
\end{array}
\end{equation}
where $G$ is a Lie group whose Lie algebra is $\g$, $\rho$ is a right representation of $G$ on $M$ and the downwards arrow is induced by a map $M\to\g$. The right $\g$-module structure on $M$ comes from its structure of a right $\mathfrak{gl}(M)$-module.  

\medskip

Note that any finite-dimensional formal group in $\LM$ comes from such a Lie group in $\LM$. 

\subsection{The universal enveloping algebra as the bimodule of 1-currents}
Finally, let us point out that the interpretation of the universal enveloping algebra as the bialgebra of distributions on a Lie group holds for the Lie groups in the Loday-Pirashvili category as defined by (\ref{LPG}). 

Indeed, let us think of a Lie group in $\LM$ as a ``generalized tangent bundle'' $\xi$ over $G$ with the fibre $M$, that is, a bundle with a two-sided action of $G$ and an anchor map to the tangent bundle. The space of sections $\Gamma(\xi^*)$ of the dual bundle $\xi^*$ can be thought of as the space of generalized differential 1-forms on $G$. Similarly, the space $D_1(\xi)$ of linear functionals on $\Gamma(\xi^*)$ that are of the form $\sum_i a_i\otimes v_i$, with $a_i\in D_0(G)$ and $v_i\in M$, can be considered as the space of generalized 1-currents supported at the unit of $G$. The anchor map sends the usual 1-forms to the generalized 1-forms and generalized 1-currents to usual 1-currents. Composing this map with $\delta: D_1(G)\to D_0(G)$ we get a map $D_1(\xi)\to D_0(G)$.

The actions of $G$ on $\xi$ give rise to a map 
$$\Gamma(\xi^*)\to \Omega^0(G)\, \widehat{\otimes}\,  \Gamma(\xi^*)+ \Gamma(\xi^*)\, \widehat{\otimes}\, \Omega^0(G). $$
By duality, $D_1(\xi)\to D_0(G)$  acquires the structure of an algebra in $\LM$. Similarly, the diagonal map $G\to G\times G$ gives rise to a map
$$\Omega^0(G)\, \widehat{\otimes}\,  \Gamma(\xi^*)+ \Gamma(\xi^*)\, \widehat{\otimes}\, \Omega^0(G)\to \Gamma(\xi^*)$$
and this gives a coalgebra structure on $D_1(\xi)\to D_0(G)$.

The primitive part of $D_1(\xi)\to D_0(G)$ is, clearly $(M\to \g)$ and it only remains to observe that the actions of $D_0(G)$ on $D_1(\xi)$ give rise to the same Lie algebra structure on $(M\to \g)$ as that coming from the right $\mathfrak{gl}(M)$ action on $M$. 

If the left action of $G$ is used to trivialize $\xi$, we can speak of generalized 1-forms with constant coefficients; these are invariant under the left action and can be identified with elements of $M^*$. The right action of $g\in G$ sends such a form $m\in M^*$ to the form $$v\to m(\mathrm{Ad}_{\rho(g)}^{-1}(v)).$$
Now, let $w\in D_0(G)$ and $v\in D_1(\xi)$ be primitive. In order to verify that $vw-wv$ coincides with the right action $[v,D\rho(w)]$ of $w$ on $v$ via $D\rho: \g\to\mathfrak{gl}(M)$ it is sufficient to check it on forms with constant coefficients. And, indeed, we have
\begin{multline*}
(vw-wv)(m)=(vw)(m)=(v\otimes w) (m(\mathrm{Ad}_{\rho(g)}^{-1}))\\ =\frac{\partial}{\partial w} m\left(\mathrm{Ad}_{\rho(g)}^{-1}(v)\right)=m\left(\frac{\partial}{\partial w} \mathrm{Ad}_{\rho(g)}^{-1}(v)\right)=m([v,D\rho(w)])=[v,D\rho(w)](m).
\end{multline*}

\begin{rem}
The solution to the integration problem given here is modelled very closely on the usual Lie theory and provides an exact analogy to the triple 
\begin{equation*}
\mbox{Lie algebras}\simeq \mbox{Irreducible cocommutative Hopf
algebras} \simeq\mbox{Formal groups},
\end{equation*}
which would be expected from any reasonable extension of Lie theory to $\LM$. The interpretation of the universal enveloping algebras in terms of 1-currents  explains in a natural way the existence of two coproducts in the enveloping dialgebras of Leibniz algebras.

Nevertheless, we should point out that the motivation behind the coquecigrue hunt is not just the quest for a good analogy, but a desire to find a homology theory for groups that would parallel the Leibniz homology for Lie algebras. The integration procedure coming from dialgebras appears to be too simple-minded to be of help in this.
\end{rem}

\begin{rem}
The same formal integration method works in other situations. In \cite{MPi} it was employed for the formal Lie theory of non-associative products. It can also be used to integrate formally Lie algebras in the category of chain complexes; the integration procedure in $\LM$ can be considered as a truncation of Lie algebra integration in that category.  
\end{rem}

\end{document}